\documentclass[12pt]{amsart}
\usepackage{amsfonts, amsmath, amssymb, amscd, latexsym,graphicx}

\newtheorem{thm}{Theorem}[section]

\newtheorem{defn}[thm]{Definition}

\def\Z{\mathbb{Z}}

\begin{document}

\title{Computational explorations in Thompson's group $F$}

\author{Jos\'e Burillo}
\address{Escola Polit{\`e}cnica Superior de Castelldefels, UPC, Avda del Canal Ol{\'\i}mpic s/n, 08860 Castelldefels, Barcelona, Spain}
\email{burillo@mat.upc.es}

\author{Sean Cleary}
\address{Department of Mathematics R8133, The City College of New York,
Convent Ave \& 138th, New York, NY 10031, USA}
\email{cleary@sci.ccny.cuny.edu}

\author{Bert Wiest}
\address{IRMAR (UMR 6625 du CNRS), Universit\'e de Rennes 1,
Campus de Beaulieu, 35042 Rennes Cedex, France}
\email{bertold.wiest@math.univ-rennes1.fr}
\thanks{We would like to thank the Centre de Recerca Matem{\`a}tica for
their hospitality and support, and the first and second authors are
grateful for support from  NSF grant \#0305545.}

\date{May 24, 2005}



\begin{abstract}
Here we describe the results of some computational explorations in
Thompson's group $F$.  We describe experiments to estimate the
cogrowth of $F$ with respect to its standard finite generating set,
designed to address the subtle and difficult question whether or
not Thompson's group is amenable.  We also describe experiments
to estimate the exponential growth rate of $F$ and the rate of escape
of symmetric random walks with respect to the standard generating set.
\end{abstract}
\maketitle

\section{Introduction}

Richard Thompson's group $F$ has attracted a great deal of interest
over the last years. The group $F$ is a finitely presented
group which arises quite naturally in different contexts, and allows
several different, but fairly simple, descriptions -- for instance
by a presentation, as a diagram group \cite{diag}, as a
group of homeomorphisms of the unit interval, as the geometry
group of associativity \cite{dehornoy}, and as the fundamental group
of a component of the loop space of the dunce hat. Cannon, Floyd and
Parry \cite{cfp} give an excellent introduction to $F$.

The interest in this group stems partly from $F$'s unusual properties,
and partly from the fact that some of the basic questions
about this group are still open, in particular those related to its
cogrowth and growth. It seems clear is that $F$
lies very close to the borderline between different regimes.

Probably the most famous open question is whether or not $F$ is amenable.
Also, it is known that $F$ has exponential growth, but the growth rate
is unknown. Similarly, the rate of escape of random walks in
$F$ is unknown.

The question of amenability is especially intriguing since $F$ is either
an example of a finitely presented non-amenable group without free
non-abelian subgroups, or an example of a finitely presented amenable but
not elementary amenable group. Though there are finitely presented
examples of groups for each of these phenomena from Grigorchuk
\cite{grig-nonelem} and Sapir and Olshanskii \cite{sapirolshan}, those
groups were constructed explicitly for those purposes, whereas $F$ is a
more ``naturally occurring'' example to consider -- so either answer
would be remarkable.

The aim of this paper is to contribute new empirical evidence to the quest to understand cogrowth,
growth, and escape rate. This evidence was obtained using large computer simulations.

The structure of this paper is as follows. In Section 2 we recall briefly
the definition and those properties of the group $F$ that will be
needed in the paper. Moreover, we give the definition of amenability
which will be used in our experiments (there are other, equivalent,
definitions which are probably more well-known). In Section 3 we
describe the algorithms used in our computations relating to amenability.
In Section 4 we present the results of our  computer experiments,
with the aim of obtaining evidence for or against the amenability of $F$.
In Section 5, we describe two computational approaches to estimate
the exponential growth rate of $F$ with respect to the standard
two-generator generating set, and in Section 6, we describe the results
of some computations to measure the average distance from the origin of
increasingly-long random walks, known as the rate of escape.


\section{Background on Thompson's group $F$ and amenability}

Richard J. Thompson's group $F$ is usually defined as the group of
piecewise-linear orientation-preserving homeomorphisms of the unit
interval, where each homeomorphism has finitely many
changes of slope (``breakpoints'') which all are dyadic integers
and and whose slopes, when defined,  are powers of 2.
$F$ admits an infinite presentation given by
$$
\left<x_1,x_2,x_3,\ldots\,|\,x_jx_i=x_ix_{j+1} \text{ if } i<j\right>
$$
which is convenient for its symmetry and simplicity, while there
is a finite presentation given by
$$
\left<x_0,x_1\,|\,[x_0x_1^{-1},x_0^{-1}x_1x_0],[x_0x_1^{-1},x_0^{-2}x_1x_0^2]\right>.
$$
Brin and Squier \cite{bs:pl}  showed that $F$ has no free
non-abelian subgroups, and thus the question of the amenability
of $F$ is potentially connected to the  conjecture of Von Neumann
that a group is amenable if and only if it had no free non-abelian
subgroups. The conjecture has since been solved negatively, but
the problem of the amenability of $F$ is of independent interest and
it has been open for at least 25 years.

The usefulness of the infinite presentation is the fact that $F$
admits a normal form based on the infinite set of generators.
The relators of the infinite presentation can be used
to reorder generators of a given word into an expression of the
following form:
$$
x_{i_1}^{r_1}x_{i_2}^{r_2}\ldots x_{i_n}^{r_n}x_{j_m}^{-s_m}\ldots
x_{j_2}^{-s_2}x_{j_1}^{-s_1}
$$
with
$$
i_1<i_2<\ldots<i_n\qquad j_1<j_2<\ldots<j_m.
$$
This normal form is unique if one requires the following extra
condition: if the generators $x_i$ and $x_i^{-1}$ both appear,
then either $x_{i+1}$ or $x_{i+1}^{-1}$ must appear as well.
Indeed, if neither $x_{i+1}$ nor $x_{i+1}^{-1}$ appeared, then
the relator could be applied so as to obtain a shorter word
representing the same element. The uniqueness of this normal form can
be used to solve the word problem in short time: given a word in the
infinite set of generators, find the normal form, which can be
done in quadratic time, and the element is the identity if and only
if the normal form is empty. This unique normal form is
most helpful when the task at hand is to decide whether two words
represent the same element of $F$. If one wishes simply to test
whether a given word represents the trivial element of $F$,
it is enough to reorder the generators, but without checking
the extra condition for uniqueness.

For an introduction to $F$ and proofs of its basic properties see
Cannon, Floyd and Parry
\cite{cfp}. Also, for an excellent introduction to amenability,
the interested reader can consult Wagon \cite{Wagon}, Chapters 10 to 12.

There are several equivalent definitions of amenability,
especially for finitely generated groups. The standard
definition is given by the existence of a finitely-additive
left-invariant probability measure on the set of subsets of $G$. If the
group is finitely generated, a celebrated characterization due to
F\o lner \cite{Folner} in terms of the existence of sets with
small boundary, has given a special interest to this concept from
the point of view of geometric group theory, making it
easier to see  that amenability is a quasi-isometry invariant.

The numerical criterion we will use extensively in this paper is due
to Kesten \cite{Kesten1,Kesten2} and it uses the concept of cogrowth.

\begin{defn} Let $G$ be a finitely generated group and let
$$
1\to K\to F_m\to G\to 1
$$
be a presentation for $G$. The \emph{cogrowth} of $G$ is the growth of the subgroup $K$ inside
$F_m$. In particular, the \emph{cogrowth function} of $G$ is
$$
g(n)=\#(B(n)\cap K),
$$
where $B(n)$ is the ball of radius $n$ in $F_m$,
and the
\emph{cogrowth rate} of $G$ is
$$
\gamma=\lim_{n\to\infty}g(n)^{1/n}.
$$
\end{defn}

Kesten's cogrowth criterion for amenability states basically that
a group is amenable when it has a large proportion of freely reduced
words, for every length $n$, representing the trivial element;
that is, when the cogrowth is large.

\begin{thm}[Kesten]\label{Kesten's criterion}
Let $G$ be a finitely generated group, and let $X$ be a finite set of generators, with cardinal
$m$. Let $\gamma$ be its cogrowth rate. Then $G$ is amenable if and only if $\gamma=2m-1$.
\end{thm}

This can also be interpreted in terms of random walks.  If the group
is nonamenable (that is, if there are very few nontrivial words representing
the trivial element of the group), then the probability of a random walk
in the group ending at 1 is small.
Since our random walks are taken to be non-reduced, we consider the
$(2m)^L$ non-reduced words of length $L$ in $m$ generators, and let
$T(L)$ be the set of these words which represent the identity in the
group $G$. Then, define
$$
p(L)=\frac{\#T(L)}{(2m)^L},
$$
that is, we define $p(L)$ to be the proportion of words which are equal
to the identity in $G$. Then, a rewriting of Kesten's criterion for
non-reduced words can be given by

\begin{thm}[Kesten] A group is amenable if and only if
$$
\limsup_{L\to\infty}p(L)^{1/L}=1
$$
\end{thm}

Roughly speaking, a group is amenable if the probability of a random
walk of length $L$ returning to $1$ decreases more slowly than
exponentially with $L$. This form of the criterion will be used in the
subsequent sections to try to study numerically the amenability of $F$.


\section{Algorithms and programs}

The direct approach at finding the numbers $p(L)$ exactly for Thompson's group $F$ fails even at
quite small values of $L$ due to the fact that the number of words grows exponentially, so the
computation times get large easily. For instance, for a length as small as 14 the number of total
words is $4^{14}$=268,435,456, out of which there are 1,988,452 representing the neutral element,
for a value $p(14)^{1/14}=0.704423677$. It would be hard to decide whether the sequence approaches
1. A number of improvements can be made to ease the calculation so it becomes more feasible to
estimate whether the sequence tends to 1.

First, we take samples of words of a given length  instead of the all
words of a given length.
The number $4^L$ grows impracticably large even for small values of $L$,
so sampling becomes a necessity. Since the number $p(L)$ is
basically a proportion (or a probability), it can be approximated by
Monte Carlo methods. One can always take a random non-reduced word in
the two generators $x_0$ and $x_1$ and check if it is the identity
by solving the word problem quickly using the normal form.
Repeating this process one can find a reasonably good approximation
of the number $p(L)$.

A further improvement can be implemented by taking only \emph{balanced}
words. We observe that, since the two relators in $G$ are commutators, a
word which represents the identity has to be \emph{balanced}: it has
to have total exponent zero in both generators $x_0$ and $x_1$.
So we consider not \emph{all} random words, but only
balanced ones. We remark that the abelianization of $F$ is $\mathbb{Z}^2$,
generated by $x_0$ and $x_1$, so being balanced is in fact equivalent to
representing the trivial element of $\Z^2=F_{\rm ab}$. Now we let $C(L)$ be
the set of balanced words among the $4^L$ non-reduced words of length $L$
in $F_2$, and define
$$
\widehat{p}(L) = \frac{\#T(L)}{\#C(L)},
$$
the proportion of words representing the identity of $F$ among balanced
words of length $L$. We have
\begin{eqnarray*}
\sqrt[L]{p(L)} & = & \sqrt[L]{\frac{\#T(L)}{4^L}}
=\sqrt[L]{\frac{\#T(L)}{\#C(L)}}\cdot \sqrt[L]{\frac{\#C(L)}{4^L}}\cr
& = & \sqrt[L]{\widehat{p}(L)} \cdot \sqrt[L]{\frac{\#C(L)}{4^L}}.
\end{eqnarray*}
Moreover, the last factor $\sqrt[L]{\frac{\#C(L)}{4^L}}$
tends to 1 as $L$ tends to infinity, because $\Z^2$ is amenable.
Thus $F$ is amenable if and only if we have
$\limsup_{L\to\infty} \widehat{p}(L)^{1/L} = 1$.

So in order to decide whether $F$ is amenable, we shall try to find good approximations of
$\widehat{p}(L)$, the proportion of words representing $1_F$ among balanced words of length $L$,
and this for values of $L$ which are as large as possible. Obviously, the algorithm for creating
random balanced words must be designed in such a way that all balanced words of length $L$ have
the same chance of appearing. The practical advantage of approximating $\widehat{p}(L)$ rather
than $p(L)$ is that $\widehat{p}(L)$ is much larger (roughly by a factor $\pi L/2$), so much
smaller sample sizes are required.

Yet another improvement, which substantially increases the efficiency
of the algorithm, can be made by using a ``divide and conquer" strategy.
The underlying observation is that if $L$ is even, then the probability
that a random word of length $L$ represents the trivial element of $F$
is equal to the probability that two random words of length $L/2$
represent the same element. Thus, the idea of the algorithm is to
create a large number $N$ of random words of length $L/2$ (in our
implementations, values for $N$
between 15,000 and 200,000 were generally used). Each of the $N$
words is immediately brought into normal form, and these normal forms
are stored. In order to decide if two words represent the same
element of $F$, we simply compare their normal forms.
Therefore we can consider all $N (N-1)/2$ unordered pairs of
words in normal form, and we count how many identical pairs we see.
This number, divided by $N(N-1)/2$, is an approximation for
the proportion $p(L)$. However, the description just provided is
an oversimplification, because as described above, we
would like to restrict our sample
to \emph{balanced} words. Here we describe the
estimation algorithm more precisely:

Each iteration of the algorithm has the following steps. In a preliminary step, we create one
random \emph{balanced} word of length $L$. Then we focus our attention on the first half (the
first $L/2$ letters) of this word and we count which element in the quotient $F_{\rm ab} = \Z^2$
this first half represents ---that is, we count the exponent sums of the letters $x_0$ and $x_1$
for the first half of the word.

In the second step, we create $N$ random words of length $L/2$ which
all represent this same element of the abelianization $F_{\rm ab}=\Z^2$,
in such a way that  all possible words of length $L/2$ with the given
$x_0$-balance and $x_1$-balance have the same chance of appearing.
As soon as it is created, each random word is transformed into normal
form, and this normal form is stored.

In the third step, we count the proportion of identical pairs among all
$N(N-1)/2$ unordered pairs of stored words in normal form.

In this way, each iteration of the algorithm gives an approximation to
the true value of $\widehat{p}(L)$. Performing a few thousand iterations,
and taking the mean of the proportions obtained in each step, one obtains
an approximation to $\widehat{p}(L)$.

The expected value for the result of this algorithm is indeed $\widehat{p}(L)$, which we interpret
as the probability that two random words of length $L/2$ represent the same element of $F$, under
the condition that they represent the same element of $F_{\rm ab}=\Z^2$.  It is immediate from the
construction of the algorithm that for any pair $(k,l)\in \Z^2$, the proportion of words
representing $(k,l)$ in $F_{\rm ab}$ among all words constructed by the algorithm is what it
should be ---namely the probability that the first half of a balanced random word of length $L$
represents $(k,l)$ in $F_{\rm ab}=\Z^2$.

Then, having fixed some pair $(k,l)$ in $\Z^2$, we restrict
our attention to those iterations of the algorithm that deal with
words with $x_0$-balance $k$ and $x_1$-balance $l$ (and length $L/2$).
We have to prove that the expected value for the proportion of identical
pairs of words in our algorithm is what it should be -- namely the
probability that a pair of random words, chosen with uniform probability
from the set pairs of words of length $L/2$ representing the element
$(k,l)$ of $F_{\rm ab}=\Z^2$, represent the same element of $F$.
That is, we have to prove that our taking words in batches of $N$ and
comparing all couples in that batch, rather than taking independent
samples of pairs of words, does not distort the result. That, however,
follows immediately from the fact that in our algorithm, all pairs of
words of length $L/2$ with $x_0$-balance $k$ and $x_1$-balance $l$,
appear on average with the same frequency (they have uniform probability).
The fact that our $N(N-1)/2$ samples are not independent
has no impact on the expected value. It does have an impact on the
variation, that is, on the size of the error bars, but even this negative
impact becomes negligible when we have, on average, less than one
identical pair per batch of $N$ words, as we typically have.

The authors have implemented the last two algorithms in computer programs
written in FORTRAN and C. These programs were run for several weeks on
the ``Wildebeest" 132-processor Beowulf cluster at the City University
of New York. The results of these implementations will be shown in the
next section.


\section{Computational results concerning amenability}

The results for the computations of trivial words for $F$ are represented
in Table \ref{f}. This table contains the following information. For
lengths $L=20, 40,\ldots, 300, 320$, it gives in the second and third
columns the sample size (the number of words that were tested) and the
number of words among them that were found to represent the trivial
element of $F$; thus the quotient of these two quantities is an
approximation of $\widehat{p}(L)$. The fourth column contains the $L$th
root of this proportion. The last column contains the $20$th root of the
quotient of the proportions obtained for length $L$ and for length $L-20$.

In order to clarify the last two columns we remark that the
sequences $\sqrt[L]{\widehat{p}(L)}$ and
$\sqrt[20]{\widehat{p}(L)/\widehat{p}(L-20)}$
have the same limits -- for instance if we had
$\widehat{p}(L)\simeq \mathrm{const}\cdot a^L$ then we would obtain
$$\lim_{L\to\infty}\sqrt[L]{\widehat{p}(L)}=
\lim_{L\to\infty}\sqrt[20]{\widehat{p}(L)/\widehat{p}(L-20)}=a$$
The difference between the two sequences is that the second one
converges much more quickly, but it is also more sensitive to
statistical errors related to insufficient sample size.

In summary, the question of amenability comes down to the question
whether the numbers in the last two columns converge to $1$, or to
a smaller number. The numbers in the second to last column converge
more slowly, but they are more reliable.

\begin{table}

\begin{center}
\begin{tabular}{|c|crcc|}
  \hline
  length $L$ & sample size & trivial & \ \ $\sqrt[L]{\widehat{p}(L)}$ &
  $\sqrt[20]{\frac{\widehat{p}(L)}{\widehat{p}(L-20)}}$ \\
  \hline
  20 & $2.000\cdot 10^7$ & 1\,364\,638 & 0.8744 &  \\
  40 & $2.000\cdot 10^7$ & 82\,922 & 0.8718 &  0.8693\\
  60 & $2.000\cdot 10^7$ & 6\,341 & 0.8744 &  0.8794\\
  80 & $2.500\cdot 10^{11}$ & 7\,255\,725 & 0.8776 &  0.8873\\
  100 & $3.125\cdot 10^{11}$ & 938\,587 & 0.8806 &  0.8928\\
  120 & $8.750\cdot 10^{12}$ & 2\,961\,321 & 0.8832 &  0.8966\\
  140 & $1.312\cdot 10^{13}$ & 551\,480 & 0.8857 &  0.9009\\
  160 & $1.238\cdot 10^{13}$ & 67\,542 & 0.8879 &  0.9030\\
  180 & $2.420\cdot 10^{13}$ & 18\,618 & 0.8900 &  0.9067\\
  200 & $1.425\cdot 10^{14}$ & 16\,040 & 0.8918 &  0.9084\\
  220 & $1.572\cdot 10^{15}$ & 26\,596 & 0.8934 &  0.9096\\
  240 & $2.063\cdot 10^{16}$ & 55\,941 & 0.8950 &  0.9125\\
  260 & $2.716\cdot 10^{16}$ & 12\,162 & 0.8964 &  0.9139\\
  280 & $7.566\cdot 10^{15}$ & 599 & 0.8976 &  0.9139\\
  300 & $1.343\cdot 10^{16}$ & 196 & 0.8993 &  0.9221\\
  320 & $5.856\cdot 10^{16}$ & 148 & 0.9003 &  0.9161\\
  \hline
\end{tabular}

\caption{Cogrowth estimates for $F$.}\label{f}
\end{center}

\end{table}

Before we can establish any conclusions, it would be interesting to
compare these results with the corresponding results for groups
which are known to be amenable or not. As test groups we will take
the free group on two generators as a nonamenable example, and the
group $\mathbb{Z}\wr\mathbb{Z}$ ($\Z$ wreath $\Z$). The latter group
is amenable since it is abelian-by-cyclic, and it appears as a subgroup
of $F$ in multiple ways \cite{diag, wreathdistort}. The group $\Z\wr\Z$
admits the presentation
$$
\left<a,t\,|\,\left[a^{t^i},a^{t^j}\right],
i,j\in\mathbb{Z}\right>,
$$
and being two-generated it appears to be a good match to compare with $F$.
The results for these two groups are in Table \ref{twogroups}.

\begin{table}

\begin{center}
\begin{tabular}{|c|crc|crc|}
\hline
&\multicolumn{3}{|c|}{$\mathbb{Z}\wr\mathbb{Z}$}&
\multicolumn{3}{|c|}{$F_2$}\\
 \cline{2-7}
  $L$ & sample & trivial & $\sqrt[L]{\widehat{p}(L)}$
  & sample & trivial & $\sqrt[L]{\widehat{p}(L)}$ \\
  \hline
20 & $2.475\cdot 10^7$ & {\tiny 1\,802\,935} & 0.8772 &
$1.000\cdot 10^7$ & {\tiny 655\,940} & 0.8727 \\
40 & $2.475\cdot 10^7$ & {\tiny 247\,710} & 0.8913 &
$1.000\cdot 10^7$ & {\tiny 30\,685} & 0.8653 \\
60 & $1.980\cdot 10^7$ & 34\,658 & 0.8996 &
$2.000\cdot 10^7$ & 2\,888 & 0.8630 \\
80 & $2.475\cdot 10^{7}$ & 9\,669 & 0.9066 &
$3.000\cdot 10^7$ & 230 & 0.8631 \\
100 & $1.980\cdot 10^{7}$ & 2\,079 & 0.9125 &
$4.000\cdot 10^8$ & 159 & 0.8630 \\
120 & $1.095\cdot 10^{8}$ & 3\,485 & 0.9173 &
$6.975\cdot 10^{11}$ & {\tiny 14\,167} & 0.8628 \\
140 & $9.950\cdot 10^{7}$ & 1\,035 & 0.9213 &
$8.000\cdot 10^{11}$ & 819 & 0.8626 \\
160 & $4.990\cdot 10^{8}$ & 1\,847 & 0.9248 &
$2.400\cdot 10^{12}$ & 136 & 0.8629 \\
180 & $2.997\cdot 10^{9}$ & 4\,141 & 0.9278 &  &  &  \\
200 & $4.740\cdot 10^{10}$ & 26\,919 & 0.9306 &  &  &  \\
220 & $8.636\cdot 10^{10}$ & 20\,625 & 0.9330 &  &  &  \\
240 & $1.859\cdot 10^{11}$ & 19\,469 & 0.9352 &  &  &  \\
260 & $4.249\cdot 10^{11}$ & 20\,112 & 0.9372 &  &  &  \\
280 & $5.734\cdot 10^{11}$ & 12\,735 & 0.9390 &  &  &  \\
300 & $5.844\cdot 10^{11}$ & 6\,256 & 0.9407 &  &  &  \\
320 & $4.050\cdot 10^{12}$ & 21\,229 & 0.9422 &  &  &  \\
 \hline
\end{tabular}
\caption{Cogrowth estimates for $\Z \wr \Z$ and the free group of rank 2.}\label{twogroups}
\end{center}
\end{table}
\begin{figure}[ht!]
  \includegraphics[width=13cm]{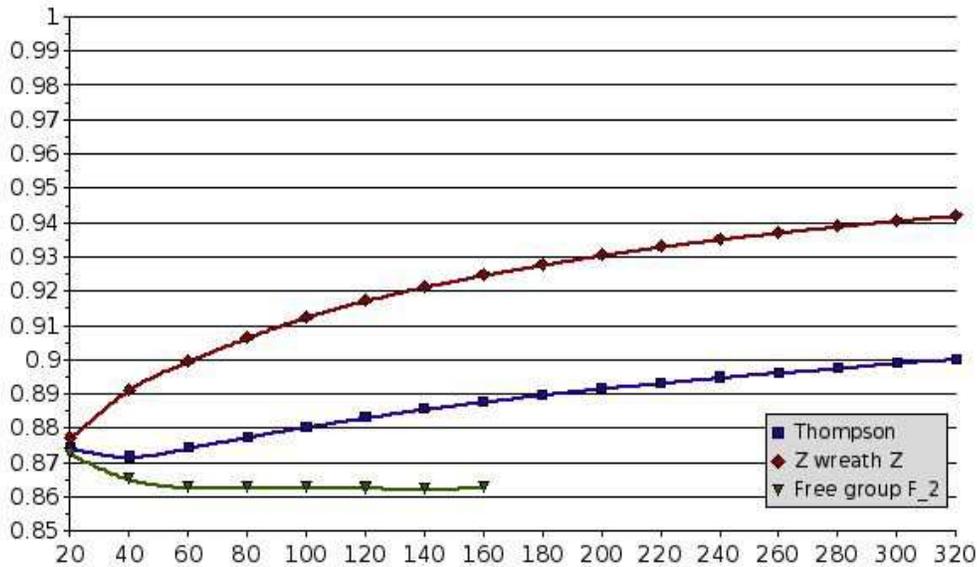}

  \caption{Comparing cogrowth estimates $\sqrt[L]{\widehat{p}(L)}$ for
  three groups.}
  \label{fig:compare}
\end{figure}

A graphical representation of  comparing these estimates of
cogrowth in the three groups $F$, $\Z \wr \Z$ and $F(2)$ is given in
Figures \ref{fig:compare}.

Do these pictures suggest that $F$ is amenable or non-amenable?
It is difficult to discern convergence to 1 or something less than 1 with
this data, and it is clear by considering other amenable groups such as
iterated  wreath products like $\Z\wr\Z\wr\Z$ that the 
convergence to 1 could be exceptionally slow.


\section{Computational results concerning the growth of $F$}

Another family of open questions about Thompson's group $F$ center
on the growth of $F$ with respect to its standard generating set
$\{x_0, x_1\}$. To study the growth of a group with respect to a
generating set, we consider $g_n$, the number of distinct elements
of $F$ of length $n$ and we form the spherical growth series,
$g(x)= \sum g_n x^n$. If we consider balls of radius $n$ and the
number of elements $b_n$ whose length is less than or equal to $n$,
we have the growth series $  b(x)= \sum b_n x^n$.
Thompson's group has exponential growth as the submonoid generated
by $x_0, x_1$ and $x_1^{-1}$ is free (see  Cannon, Floyd and Parry \cite{cfp}).
Burillo \cite{josegrowth} computed the exact growth function for positive
words in $F$ with respect to the standard two generator generating set
$\{x_0,x_1\}$ which gives a lower bound for the growth rate of words
in the full group as the largest root of $x^3-2x^2-x+1$, which is about 2.24698.
Guba \cite{gubagrowth} used the normal forms for elements of
$F$ developed by Guba and Sapir \cite{gubasapirdehn} to sharpen the lower
bound of the growth function to $\frac12 (3+\sqrt 5)$ which is
about 2.61803. Guba  conjectures that  2.7956043 is an upper bound
by considering the ratio of the ninth and eighth terms in the spherical
growth series of $F$.  But the exact growth function of $F$ remains
unknown -- it is not even known if the growth function is rational, though
Cleary, Elder and Taback \cite{lamplang} show that there are infinitely
many cone types, which may be evidence that the growth of the full
language of geodesics is not rational.

Here, we use a computational approach to estimate the growth function
of $F$. We use two methods both based upon taking random samples of words
via random walks.
Both of these methods estimate the number of words in successive
$n$-spheres of $F$. For the first method, we take an element of length $n$
and consider its ``inward'' and ``outward'' valence in the Cayley graph.
Since the relators of $F$ with respect to the standard finite presentation
are all of even length, application of a generator $x$ to an element $w$
of $F$ will either increase or reduce the length by 1. The
{\em inward valence} of $w$ is the number of generators which reduce
the word length and the {\em outward valence} of $w$ is the number of
generators which increase word length.  If the length of $w$ is $n$, then
the outward valence gives the number of words adjacent to $w$ which lie
on the $n+1$ sphere.  By taking an average of the outward valence of a
large number of elements in the
$n$ sphere, we can estimate the ratio of the number of elements in
the $n+1$ sphere to the number of elements in the $n$ sphere. Thus we
can estimate the rate of growth, as the limit of these ratios (for
$n\to \infty$) will be the exponential growth rate for the group.

For the second method, we consider a variation of this approach where
instead of looking at the words at distance 1 from $w$, we look at the
words at distance 2 from $w$ and see how many of those words lie in
the $n+2$ sphere. This gives an estimate of the ratio of the number of
elements in the $n+2$ sphere to the number of elements in the $n$ sphere,
and in the limit, we expect the square root of these
ratios to approach the exponential growth rate for the group.

We expect both methods to yield overestimates of the true growth rate,
but the error should be larger for the first method than for the second
one. The raw outward valence method is expected to overestimate because it
may count elements in the $n+1$ sphere which are adjacent to more than one
element in the $n$ sphere multiple times.  An extreme example of this are
``dead-end'' elements in $F$, characterized by Cleary and Taback
\cite{ctcomb}.  These dead-end elements have the property that right
multiplication by any generator reduces word length.  The ``outward
valence'' method includes these dead-end elements in the count of
growth -- if the randomly selected element in the $n$ sphere is one of
the 4 elements in the $n$ sphere which is adjacent to a particular
dead-end element in the $n+1$ sphere, it will contribute to the average
outward valence at least 1. For the distance two method, however, such
elements will not contribute to the growth as there will be no words
adjacent to the dead-end element which lie in the $n+2$ ball.

To compute the length of an element of $F$, we use Fordham's method
\cite{blakegd} for measuring word length of elements of $F$ with respect
to $\{x_0, x_1\}$.  This remarkable method amounts to building the reduced
tree pair diagram associated to an element of $F$, classifying each
internal node of the trees diagram into one of seven possible types, and
then pairing the nodes and summing a weight function of those node type
pairs to get the exact length of the element.

\begin{figure}[htb]
  \includegraphics[width=12cm]{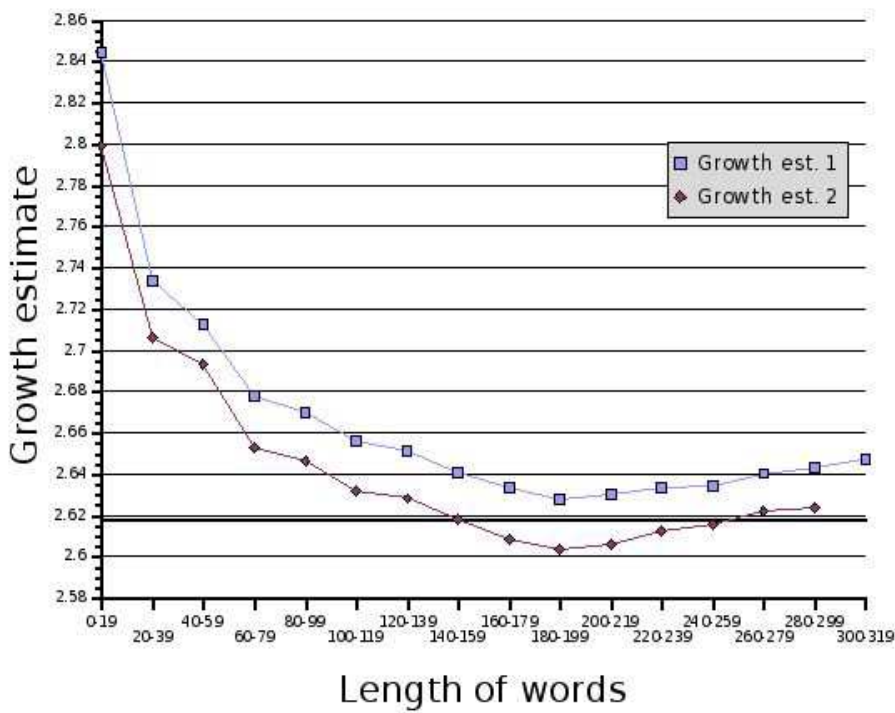}
  \caption{Estimates for the exponential growth rate from the data in
  Table \ref{growthtable}}
  \label{fig:growth}
\end{figure}

We note that selecting a random element of
the $n$ sphere for a predetermined value of $n$ is not
feasible given current understanding of the metric balls in $F$ -- we do
not even know the number of such elements, as in fact that is
what we are trying to estimate.  So we construct elements by taking
random walks in the group with respect to the standard generating set of
a predetermined length $n$, and then measure the length $l$ of the element
obtained.  We then compute its outward valence by measuring the lengths of
elements adjacent to it in the Cayley graph and we also count the number
of elements at distance two from it which lie in the $l+2$ sphere. Thus,
we obtain simultaneously estimates of outward valence for elements in a
range of balls.  Furthermore, we can record the length $l$ of a word
obtained by a random walk of length $n$ and use that to estimate crudely
the rate of escape of a random walk in $F$, as described in the next
section. The results of the computations concerning growth are presented
in Table \ref{growthtable} and Figure \ref{fig:growth}.

\begin{table}
\begin{center}
\begin{tabular}{|c|c|c|c|c|}
\hline
 Lengths & Words & \parbox{1in}{Average outward valence}
  & \parbox{1in}{Average num. at  dist. 2} & \parbox{1in}{Growth estimate from dist 2} \\
  \hline
 0 - 19   & 5723    &  2.8440    &  7.8363 & 2.7993 \\
20 - 39   & 629964  &  2.7334    &  7.3239 & 2.7063 \\
40 - 59   & 1017998 &  2.7128    & 7.2521  & 2.6930 \\
60 - 79   & 602694  &  2.6781    & 7.0389  & 2.6531 \\
80 - 99   & 612613  &  2.6698    & 7.0041  & 2.6465 \\
100 - 119 & 514665  &  2.6564    & 6.9256  & 2.6317 \\
120 - 139 & 392069  &  2.6512    & 6.9074  & 2.6282 \\
140 - 159 & 272564  &  2.6407    & 6.8529  & 2.6178 \\
160 - 179 & 234893  &  2.6331    & 6.8057  & 2.6088 \\
180 - 199 & 281806  &  2.6275    & 6.7779  & 2.6034 \\
200 - 219 & 283764  &  2.6299    & 6.7897  & 2.6057 \\
220 - 239 & 164359  &  2.6336    & 6.8234  & 2.6122 \\
240 - 259 & 48750   &  2.6341    & 6.8431  & 2.6159 \\
260 - 279 & 7326    &  2.6403    & 6.8756  & 2.6221 \\
280 - 299 & 521     &  2.6430    & 6.8829  & 2.6235 \\
300 - 319 & 17      &  2.6470    & 6.8235  & 2.6122 \\

  \hline
\end{tabular}
\end{center}
\caption{Average outward valence of words arising from random walks.}
\label{growthtable}
\end{table}

As we can see from the data, and as expected, the estimates using the
distance two method are lower than the estimate from the outward
valence method. Moreover, for the first experiment, the values lie
between the proven lower bound of 2.618\ldots and the conjectured upper
bound of 2.763\ldots, for words of length 20 and more.
However, other aspects of the computational results are more surprising.
Both functions appear to have a minimum at length about 190. Moreover,
for the second experiment, the values obtained lie below the
proven lower bound for words of length between 140 and 260, and
lie in the expected range before and after that.   This data suggests
that the rate of growth is close to the proven lower bound or that
random walks are not an unbiased method for estimating growth
by average outward valence. 
 Of course, since we do not know
the growth function, it is difficult to effectively pick a random element,
so perhaps random walks tends to bias
toward those which have lower outward valence than is representative.
The role of ``dead-end'' elements of outward valence 0 may play a role
in this bias and we describe estimates of densities of dead-end elements in the next section.
It may be that random walks get stuck near dead-end elements and other low outward
valence items and thus random walks may select these elements at a greater
proportion than uniform.

Finally, we mention that we have also computed first twelve terms of
the exact spherical growth function of $F$ to obtain:
\begin{eqnarray*}
g(x) & = & 1+ 4x+ 12x^2 +36x^3 +108x^4 + 314x^5+ 906 x^6 + 2576 x^7 +\\
 & & \ \ \ \ +  7280 x^8  +  20352x^9+ 56664x^{10} +  156570 x^{11} + \ldots
\end{eqnarray*}

Guba \cite{gubagrowth} had already calculated the first ten terms of this sequence and
noticed that the ratios of successive terms of this series
appear to decrease and form a natural conjectural upper bound to the growth function. The two
additional successive quotients arising from our additional terms continue the decreasing pattern
and are $2.7841981\ldots$ and $2.7631300\ldots$ and lie well above the experimental
estimates of growth described above.


\section{Rate of escape of random walks and dead-ends in $F$}

Here we note that as a side effect of the computations described in
the previous section to estimate growth, we obtain two pieces of data
which are interesting in their own right.

First, since the random elements used to estimate growth are constructed
by random walks and we measure their exact lengths using Fordham's method,
we are able to see how quickly these random walks leave the origin.
Since these are symmetric random walks, there is of course the possibility
of backtracking to get non-freely reduced words, so we do not expect a
random walk of length 100 to actually reach the sphere of radius 100 with
non-negligible probability. Our estimates of the rate of escape of random
walks of lengths 100 to 1000 are shown in Table \ref{tab:escape} and the
rate of escape seems to be decreasing in this range.

\begin{table}[htb]
\begin{center}
\begin{tabular}{|c|c|c|c|c|}
\hline
\parbox{.7in}{ Length of random walk } & \parbox{.7in}{Number of walks}  & \parbox{.6in}{Average length}  &  \parbox{.7in}{Standard deviation} & \parbox{.5in}{Rate of escape}\\
\hline
  100& 4764000 & 41.18 & 8.34 & 0.4118\\
  200& 3242898 & 76.01 & 12.33 & 0.3800 \\
300  & 2700000 & 109.3 & 15.51& 0.3545 \\
400 &1500000 & 141.8 & 18.33 & 0.3544 \\
500 &600000 & 173.8 & 20.82 &0.3476 \\
600 & 1500000 & 205.3 & 23.08 & 0.3421\\
700 & 900000 & 236.5 & 25.14 &0.3379 \\
800  & 900000 & 267.6 & 27.14& 0.3345 \\
900 & 300000 & 298.5 & 29.02 &0.3316 \\
1000 & 300000 & 329.0 & 30.86 &0.3290 \\
\hline
\end{tabular}
\end{center}
\caption{Distance from origin (word length) as a function of random walk length}
\label{tab:escape}
\end{table}

Second, since we compute the outward valence of words to estimate the
growth, we can look for words of outward valence zero- these are exactly
the ``dead-end'' elements discovered by Fordham \cite{blake:diss} and
characterized by Cleary and Taback \cite{ctcomb}. Though dead-end elements
can occur in any group (with respect to generating sets contrived for that
purpose) groups with dead-end elements with respect to natural generating
sets are much less common. Geodesic rays from the identity towards infinity
cannot pass through dead-end elements, and thus the existence of many
dead-end elements tends to reduce the growth of the group.
Table \ref{tab:dead} shows the observed incidence of dead ends during the
course of the growth estimation calculations in Section 5.  We see that
there are significant numbers of dead ends but that the fraction decreases
as the lengths of elements increases.

\begin{table}[htb]
\begin{center}
\begin{tabular}{|c|c|c|c|}
\hline
 Range of lengths   &  Number of words  & \parbox{.8in}{Number of dead-ends}&  Fraction \\
  \hline
0 - 39 & 634927 & 665 & 0.001047  \\
40 - 79 & 1620692 & 1386 & 0.0008552  \\
80 - 119 & 1127278 & 625 & 0.0005544  \\
120 - 159 & 665245 & 239 & 0.0003593  \\
160 - 199 & 561502 & 149 & 0.0002654  \\
200 - 239 & 825785 & 162 & 0.0001962  \\
240 - 279 & 689500 & 114 & 0.0001653  \\
280 - 319 & 393643 & 39 &   0.00009907  \\
320 - 359 & 128254 & 11 &   0.00008577  \\
360 - 399 & 20926 & 1 &        0.00004779  \\
400 - 439 & 1193 & 0 & 0  \\
440 - 479 & 21 & 0 & 0  \\
\hline
\end{tabular}
\end{center}
\caption{Fractions of dead-ends observed during random walks as a function
of resulting word length.}
\label{tab:dead}
\end{table}

\bibliographystyle{plain}
\def\cprime{$'$}

\end{document}